\documentclass[12pt]{article}
\usepackage{textcomp}
\usepackage{graphicx}           
\usepackage{epstopdf}
\usepackage{color}              
\usepackage{xcolor,graphicx}
\usepackage{geometry}
\usepackage{float}
\usepackage{indentfirst}        
\usepackage{amsmath,amssymb,bm} 
\usepackage{amsthm}
\usepackage{cases}              
\usepackage{setspace}
\usepackage{hyperref}
\usepackage{titlesec}
\usepackage[all]{xy}            
\usepackage{tikz}
\usepackage{multirow}
\usepackage{graphicx} 
\usepackage{epstopdf}
\usepackage{caption}
\usepackage{diagbox}
\usepackage{mathrsfs}
\usepackage{tikz}
\usepackage{dsfont}
\usepackage{graphics}
\usepackage{amsmath,amssymb,amsfonts,amsthm,bm,mathrsfs}
\usepackage{tcolorbox}
\usepackage{fancyhdr,fancyvrb}
\usepackage{xcolor,graphicx,geometry,float,indentfirst,cases,setspace,titlesec}
\usepackage{diagbox}

\newtheorem{thm}{Theorem}[section]

\newtheorem{lem}[thm]{Lemma}

\newtheorem{coro}[thm]{Corollary}

\newtheorem{pf}{proof}[section]
\theoremstyle{remark}

\begingroup
\makeatletter \@addtoreset{equation}{section} \makeatother
\makeindex 
\def\qed{\hfill \rule{4pt}{7pt}}
\def\pf{\vskip 0.2cm {\noindent \bf Proof.}\quad}

\begin{document}

\begin{center}
	
	{\Large \bf Unimodality of $k$-Regular  Partitions\\[5pt] into Distinct Parts with Bounded Largest Part }
	
\end{center}

\begin{center}
	{ Janet J.W. Dong}$^{1}$ and {Kathy Q. Ji}$^{2}$  \vskip 2mm

	$^{1,2}$ Center for Applied Mathematics,  Tianjin University, Tianjin 300072, P.R. China\\[6pt]
	\vskip 2mm
	
	$^1$dongjinwei@tju.edu.cn and $^2$kathyji@tju.edu.cn
\end{center}

\vskip 6mm \noindent {\bf Abstract.}  A $k$-regular partition into distinct parts is a partition into distinct parts  with no part divisible by $k$.  In this paper, we provide  a general method to establish the unimodality of  $k$-regular  partition  into distinct parts where the largest part is at most $km+k-1$. Let $d_{k,m}(n)$ denote the number of   $k$-regular  partition of $n$  into distinct parts where the largest part is at most $km+k-1$.   In line with this method, we show that   $d_{4,m}(n)\geq d_{4,m}(n-1)$ for $m\geq 0$, $1\leq n\leq 3(m+1)^2$ and $n\neq 4$ and $d_{8,m}(n)\geq d_{8,m}(n-1)$ for $m\geq 2$ and $1\leq n\leq 14(m+1)^2$.  When $5\leq k\leq 10$ and $k\neq 8$, we show that $d_{k,m}(n)\geq d_{k,m}(n-1)$ for $m\geq 0$ and $1\leq n\leq \left\lfloor\frac{k(k-1)(m+1)^2}{4}\right\rfloor$.

\noindent
{\bf Keywords:} Unimodal, symmetry,  integer partitions, $k$-regular partitions, analytical method

\noindent
{\bf AMS Classification:} 05A17, 05A20, 11P80, 41A10, 41A58

\vskip 6mm

\section{Introduction}

The main theme of  this paper is to investigate unimodality of $k$-regular partition into distinct parts where the largest part is at most $km+k-1$.
A $k$-regular partition into distinct parts is a partition into distinct parts  with no part divisible by $k$.
For example, below are the $4$-regular partitions of $10$ into distinct parts,
\[(10), (9,1),\ (7,3),\ (7,2,1),\ (6,3,1),\ (5,3,2).\]
Let $d_{k,m}(n)$ denote the number of $k$-regular partition into distinct parts where   the largest part is at most $km+k-1$.
From the example above, we see that  $d_{4,1}(10)=4$. By definition, it is easy to see that the generating function of  $d_{k,m}(n)$ is given by
\begin{align}\label{eq-main}
	D_{k,m}(q):=\sum_{m=0}^{N(k,m)}d_{k,m}(n)q^n=\prod_{j=0}^m \left(1+q^{kj+1}\right)\left(1+q^{kj+2}\right)\cdots\left(1+q^{kj+k-1}\right),
\end{align}
where
\[N(k,m)= \frac{k(k-1)(m+1)^2}{2}.\]

Recall that a polynomial $a_0+a_1q+\cdots+a_Nq^N$
with integer coefficients is called  unimodal if for some $0\leq j\leq N$,
\[a_0\leq a_1\leq \cdots \leq a_j\geq
a_{j+1}\geq \cdots \geq a_N,\]
and is called symmetric if for all $0\leq j\leq N$, $a_j=a_{N-j},$ see \cite[p. 124, Ex. 50]{Stanley-1997}. It is well-known that the  Gaussian polynomials
\[{n\brack k}=\frac{(1-q^n)(1-q^{n-1})\cdots(1-q^{n-k+1})}{(1-q)(1-q^2)\cdots(1-q^k)}\]
are symmetric and unimodal, as  conjectured by Cayley \cite{Cayley-1856} in 1856  and confirmed by Sylvester \cite{Sylvester-1878} in 1878 based on semi-invariants of binary forms. For more information, we refer to  \cite{Chen-Jia-2022, O'Hara-1990,   Pak-Panova-2013, Proctor-1982}. Since then, the unimodality of polynomials (or combinatorial sequences) has drawn great attention  in recent decades. In particular, the unimodality of several  special $k$-regular partitions have been investigated by several authors.  For example, the polynomials
\begin{equation}\label{eq2-Odlyzko-Richmond}
	(1+q)(1+q^2)\cdots(1+q^m)
\end{equation}
are proved to be symmetric and unimodal for  $m\geq 1$.  The first   proof of the unimodality of the polynomials \eqref{eq2-Odlyzko-Richmond}  was given  by Hughes \cite{Hughes-1977} resorting to  Lie algebra results. Stanley \cite{Stanley-1982}  provided an alternative proof by using the Hard Lefschetz Theorem.  Stanley \cite{Stanley-1980} also established the general result of this type  based on a   result of Dynkin  \cite{Dynkin-1950}. An analytic proof of  the unimodality of the polynomials \eqref{eq2-Odlyzko-Richmond} was attributed to Odlyzko and Richmond \cite{Odlyzko-Richmond-1982} by  extending the argument of van Lint \cite{Lint-1967} and  Entringer \cite{Entringer-1968}.

Stanley \cite{Stanley-1982} conjectured
the polynomials
\begin{equation}\label{eq2-Almkvist}
	(1+q)(1+q^3)\cdots(1+q^{2m+1})
\end{equation}
are unimodal for $m\geq 26$, except at the coefficient of $q^2$ and $q^{(m+1)^2-2}$, which has been proved by Almkvist \cite{Almkvist-1985} via
refining the method of Odlyzko and Richmond \cite{Odlyzko-Richmond-1982}.  Pak and Panova \cite{Pak-Panova-2014} showed that the polynomials \eqref{eq2-Almkvist} are strict unimodal  by interpreting
the differences between numbers of certain partitions as Kronecker coefficients of representations of $S_n$.
By refining the method of Odlyzko and Richmond \cite{Odlyzko-Richmond-1982}, we show that the polynomials
\begin{equation}\label{dong-ji}
	\prod\limits_{j=0}^{m}(1+q^{3j+1})(1+q^{3j+2})
\end{equation}
are symmetric and unimodal for $m\geq 0$, see  \cite{Dong-Ji-2023}.

In this paper, we aim to establish the symmetry and   unimodality of $D_{k,m}(q)$ for $k\geq 4$. It should be noted that the polynomial \eqref{eq2-Odlyzko-Richmond}  is associated with $D_{1,m}(q)$, while the polynomial \eqref{eq2-Almkvist}  is associated with $D_{2,m}(q)$.   When $k=3$, $D_{k,m}(q)$ reduces to the polynomial \eqref{dong-ji}.

One main result of this paper is to  show that $D_{4,m}(q)$ is almost unimodal.

\begin{thm}\label{thm-k4}
	The polynomials
	\begin{equation}\label{dong-ji-4}
		\prod\limits_{j=0}^{m}(1+q^{4j+1})(1+q^{4j+2})(1+q^{4j+3})
	\end{equation}
	are  unimodal for $m\geq 0$, except at the coefficient of $q^4$ and $q^{6(m+1)^2-4}$.
\end{thm}

We also provide an effective way to establish the unimodality of $D_{k,m}(q)$ for $k\geq 5$.

\begin{thm}\label{main-thm}
	For $k\geq5$, if there exists $m_0\geq 0$ such that
	$D_{k,m_0}(q)$ is unimodal and  for $m_0< m< 8k^{\frac{3}{2}}$ and $\left\lceil\frac{k(k-1)m^2}{4}\right\rceil\leq n\leq\left\lfloor\frac{k(k-1)(m+1)^2}{4}\right\rfloor$,
	\begin{equation}\label{eq-4.2}
		d_{k,m}(n)\geq d_{k,m}(n-1),
	\end{equation}
	then $D_{k,m}(q)$ is unimodal for $m\geq m_0$.
\end{thm}

By utilizing Theorem  \ref{main-thm} and conducting tests with Maple, we obtain the following two consequences.

\begin{coro} \label{coro-k5}
	When $5\leq k\leq 10$ and $k\neq 8$, the polynomials
	\[  \prod_{j=0}^m \left(1+q^{kj+1}\right)\left(1+q^{kj+2}\right)\cdots\left(1+q^{kj+k-1}\right)
	\]
	are unimodal for $m\geq 0$.
\end{coro}

\begin{coro} \label{coro-k8}
	The polynomials
	\[  \prod_{j=0}^m \left(1+q^{8j+1}\right)\left(1+q^{8j+2}\right)\cdots\left(1+q^{8j+7}\right)
	\]
	are   unimodal for $m\geq 2$.
\end{coro}

\section{A Key Lemma}

This section is devoted to the proof of the following lemma. It turns out that this lemma   figures prominently in the proofs of Theorem \ref{thm-k4} and Theorem \ref{main-thm}.

\begin{lem}\label{main-lemma}
	If $k\geq4$, $m\geq 8k^{\frac{3}{2}}$ and $\frac{k(k-1)m^2}{4}\leq n\leq\frac{k(k-1)(m+1)^2}{4}$, then
	\begin{equation}\label{eq-4.2}
		d_{k,m}(n)> d_{k,m}(n-1).
	\end{equation}
\end{lem}

Before demonstrating  Lemma \ref{main-lemma},   we collect several identities and inequalities which will be useful in its proof.
\begin{align}
	&\mathrm{e}^{ix}=\cos(x)+i \sin(x),\label{eq-2.1}\\[5pt]
	&\cos (2x)=2\cos^2(x)-1,\label{eq-2.2}\\[5pt]
	&\sin (2x)=2\sin (x) \cos (x), \label{eq-2.2-3}\\[5pt]
	&2\sin(\alpha)\cos(\beta)=\sin(\alpha+\beta)+\sin(\alpha-\beta),\label{eq-2.3}\\[5pt]
	&\sin (x)\geq x\mathrm{e}^{- x^2/3}\quad\text{for} \ 0\leq x\leq2,\label{eq-2.4}\\[5pt]
	&\cos (x)\geq\mathrm{e}^{-\gamma x^2}\quad \text{for} \ |x|\leq 1,\,(\gamma =-\log \cos(1)=0.615626\ldots),\label{eq-2.5}\\[5pt]
	&x-\frac{x^3}{6} \leq \sin (x)\leq x\quad \text{for} \ x\geq 0,\label{eq-2.6}\\[5pt]
	&\cos (x)\leq  \mathrm{e}^{-x^2/2}\quad\text{for} \ |x|\leq \frac{\pi}{2},\label{eq-2.7-1}\\[5pt]
	&|\cos (x)|\leq \exp\left(-\frac{1}{2}\sin^2(x)-\frac{1}{4}\sin^4 (x)\right),\label{eq-2.7}\\[5pt]
	&\left|\frac{\sin(nx)}{\sin(x)}\right|\leq n \quad \text{for} \ x\neq i\pi,\,i=0,1,2,\ldots,\label{eq-2.7-2}\\[5pt]
	&\sum_{k=1}^n\sin^2(kx)=\frac{n}{2}-\frac{\sin((2n+1)x)}{4\sin (x)}+\frac{1}{4} \quad \text{for} \ x\neq i\pi,\,i=0,1,2,\ldots,\label{eq-2.8}\\[5pt]
	&\sum_{k=1}^n\sin^4(kx)=\frac{3n}{8}-\frac{\sin((2n+1)x)}{4\sin (x)}+\frac{\sin((2n+1)2x)}{16\sin(2x)}+\frac{3}{16} \nonumber\\[5pt] &~
	~~~~~~~~~~~~~~~~~
	~~~~~~~~~~~~~~~~~~
	~~~~~~~~~~~~~~~~~~\text{for} \ x\neq \frac{i\pi}{2},\,i=0,1,2,\ldots.\label{eq-2.9}
\end{align}
The identity \eqref{eq-2.1} is Euler's identity, see \cite[p. 4]{Stein-Shakarchi-2003}. The formulas \eqref{eq-2.2}--\eqref{eq-2.3}  of trigonometric functions can be found in \cite[Chap. 8] {Burdette-1973}. The inequalities \eqref{eq-2.4}--\eqref{eq-2.7-2} are due to Odlyzko and Richmond \cite[p. 81]{Odlyzko-Richmond-1982}.
The identities \eqref{eq-2.8} and \eqref{eq-2.9} have been proved   in \cite{Dong-Ji-2023}.

We are now in a position to prove Lemma \ref{main-lemma} by considering $d_{k,m}(n)$ as the Fourier coefficients of $D_{k,m}(q)$ and proceeding to estimate its integral.

\noindent{\it Proof of  Lemma \ref{main-lemma}: }
Putting $q=\mathrm{e}^{2i\theta}$ in \eqref{eq-main}, we get
\begin{align}
	D_{k,m}(\mathrm{e}^{2i\theta})
	&=\prod_{j=0}^{m}(1+(\mathrm{e}^{2i\theta})^{jk+1})(1+(\mathrm{e}^{2i\theta})^{jk+2})\cdots(1+(\mathrm{e}^{2i\theta})^{jk+k-1})\nonumber\\[5pt]
	&\overset{\eqref{eq-2.1}}{=}\prod_{j=0}^{m}\prod_{l=1}^{k-1}\left(1+\cos\left(2(jk+l)\theta\right)+i\sin(2(jk+l)\theta)\right)\nonumber\\[5pt]
	&\overset{\eqref{eq-2.2}\&\eqref{eq-2.2-3}}{=}\prod_{j=0}^{m}\prod_{l=1}^{k-1}\left(2\cos^2((jk+l)\theta)+2i\sin((jk+l)\theta)\cos((jk+l)\theta)\right)\nonumber\\[5pt]
	&\overset{\eqref{eq-2.1}}=\prod_{j=0}^{m}\prod_{l=1}^{k-1}2\cos((jk+l)\theta)\exp(i(jk+l)\theta)\nonumber\\[5pt]
	&=2^{(k-1)(m+1)}\exp(iN(k,m)\theta)\prod_{j=0}^{m}\prod_{l=1}^{k-1}\cos((jk+l)\theta).
	\label{eq-4.3}
\end{align}
Using Taylor's theorem \cite[pp. 47--49]{Stein-Shakarchi-2003}, we derive that
\begin{align}
	d_{k,m}(n)
	&=\frac{1}{2\pi i}\int_{-\frac{\pi}{2}}^{\frac{\pi}{2}}\frac{D_{k,m}\left(\mathrm{e}^{2i\theta}\right)}{\left(\mathrm{e}^{2i\theta}\right)^{n+1}}\mathrm{d}\left(\mathrm{e}^{2i\theta}\right)\nonumber\\[5pt]
	&=\frac{1}{\pi}\int_{-\frac{\pi}{2}}^{\frac{\pi}{2}}D_{k,m}\left(\mathrm{e}^{2i\theta}\right)\mathrm{e}^{-2in\theta}\mathrm{d}\theta\nonumber\\[5pt]
	&\overset{\eqref{eq-4.3}}=\frac{2^{(k-1)(m+1)}}{\pi}\int_{-\frac{\pi}{2}}^{\frac{\pi}{2}}\exp(i(N(k,m)-2n)\theta)\prod_{j=0}^{m}\prod_{l=1}^{k-1}\cos((jk+l)\theta)\mathrm{d}\theta\nonumber\\[5pt]
	&\overset{\eqref{eq-2.1}}{=}\frac{2^{(k-1)(m+1)}}{\pi}\int_{-\frac{\pi}{2}}^{\frac{\pi}{2}}\left(\cos((N(k,m)-2n)\theta)+i\sin((N(k,m)-2n)\theta)\right)\nonumber\\[5pt]
	&\hskip 4cm \times \prod_{j=0}^{m}\prod_{l=1}^{k-1}\cos((jk+l)\theta)\mathrm{d}\theta.\nonumber
\end{align}	
Observe that
\[\int_{-\frac{\pi}{2}}^{\frac{\pi}{2}}\sin((N(k,m)-2n)\theta)\prod_{j=0}^{m}\prod_{l=1}^{k-1}\cos((jk+l)\theta)\mathrm{d}\theta=0,\]
so we conclude that
\[d_{k,m}(n)=\frac{2^{(k-1)(m+1)+1}}{\pi}\int_{0}^{\frac{\pi}{2}}\cos((N(k,m)-2n)\theta)\prod_{j=0}^{m}\prod_{l=1}^{k-1}\cos((jk+l)\theta)\mathrm{d}\theta.\]
To show that $d_{k,m}(n)$ increases with $n$, we take the derivative with respect to $n$,
\[\frac{\partial}{\partial n}d_{k,m}(n)=\frac{2^{(k-1)(m+1)+2}}{\pi}\int_0^{\frac{\pi}{2}}\theta\sin\left((N(k,m)-2n)\theta\right)\prod_{j=0}^{m}\prod_{l=1}^{k-1}\cos((jk+l)\theta)\mathrm{d}\theta.\]
Let $N(k,m)-2n=\mu$, and let
\[I_{k,m}(\mu)=\int_0^{\frac{\pi}{2}}\theta\sin\left(\mu\theta\right)\prod_{j=0}^{m}\prod_{l=1}^{k-1}\cos((jk+l)\theta)\mathrm{d}\theta.\]
Thus it suffices to show that
\begin{equation}\label{eq-4.4}
	I_{k,m}(\mu)> 0 \ \text{for} \ k\geq 4,\ m\geq 8k^{\frac{3}{2}}\ \text{and} \ 0<\mu\leq\frac{k(k-1)(2m+1)}{2}.
\end{equation}
We will separate the integral $I_{k,m}(\mu)$ into three parts,
\begin{align}
	I_{k,m}(\mu)
	&=\left\{\int_0^{\frac{2\pi}{k(k-1)(2m+1)}}+\int_{\frac{2\pi}{k(k-1)(2m+1)}}^{\frac{\pi}{2km+2(k-1)}}+\int_{\frac{\pi}{2km+2(k-1)}}^{\frac{\pi}{2}}\right\}\theta\sin\left(\mu\theta\right)\prod_{j=0}^{m}\prod_{l=1}^{k-1}\cos((rk+j)\theta)\mathrm{d}\theta\nonumber\\[5pt]
	&=I_{k,m}^{(1)}(\mu)+I_{k,m}^{(2)}(\mu)+I_{k,m}^{(3)}(\mu),\nonumber
\end{align}	
and  show that when $k\geq 4,\ m\geq 8k^{\frac{3}{2}}$ and $ 0<\mu\leq\frac{k(k-1)(2m+1)}{2}$,
\begin{equation}\label{eq-4.5}
	I_{k,m}^{(1)}(\mu) >
	\left|I_{k,m}^{(2)}(\mu)\right|+\left|I_{k,m}^{(3)}(\mu)\right|,
\end{equation}
from which, it is immediate that \eqref{eq-4.4} is valid.

We first estimate the value of $I_{k,m}^{(1)}(\mu)$. Recall that
\begin{equation}\label{defi-I1} I_{k,m}^{(1)}(\mu)=\int_0^{\frac{2\pi}{k(k-1)(2m+1)}}\theta\sin\left(\mu\theta\right)\prod_{j=0}^{m}\prod_{l=1}^{k-1}\cos((jk+l)\theta)\mathrm{d}\theta.
\end{equation}
When $0\leq \theta\leq \frac{4}{k(k-1)(2m+1)}$, we see that $0\leq \mu \theta \leq 2 \ \text{and}\ 0\leq(jk+l)\theta\leq 1 \ \text{for}  \ 0\leq j\leq m \ \text{and} \ 1\leq l\leq k-1$. Using \eqref{eq-2.4} and \eqref{eq-2.5}, we deduce that  \[\sin\left(\mu\theta\right)\geq \mu \theta \exp\left(-\frac{\mu^2\theta^2}{3}\right) \ \text{and} \ \cos((jk+l)\theta)\geq \exp\left(-\gamma((jk+l)\theta)^2)\right).\]
Hence
\begin{align*}
	&\theta\sin\left(\mu\theta\right)\prod_{j=0}^{m}\prod_{l=1}^{k-1}\cos((jk+l)\theta)\nonumber\\[5pt]
	&\quad\geq \mu\theta^2\exp\left(-\frac{\mu^2\theta^2}{3}\right)
	\exp\left(-\gamma\theta^2\sum_{j=0}^m\sum_{l=1}^{k-1}(jk+l)^2\right)\nonumber\\[5pt]
	&\quad\geq \mu\theta^2\exp\left(-\frac{k^2(k-1)^2(m+\frac{1}{2})^2\theta^2}{3}\right)\nonumber\\[5pt]
	&\qquad\times\exp\left(-\gamma\theta^2k(k-1)\left(\frac{km^3}{3}+km^2+\frac{(6k-1)m}{6}+\frac{2k-1}{6}\right)\right).
\end{align*}	
Put
\begin{equation*}
	c_k(m)=k^2(k-1)^2\left(\frac{1}{3m}+\frac{1}{3m^2}+\frac{1}{12m^3}\right)+\gamma k(k-1)\left(\frac{k}{3}+\frac{k}{m}+\frac{6k-1}{6m^2}+\frac{2k-1}{6m^3}\right).
	\label{eq-c_r(n)}
\end{equation*}
When $k\geq4$ and $m\geq 8k^{\frac{3}{2}}$, we find that
\begin{align*}\label{defi-cr}
	c_k(m)&\leq c_k\left(8k^{\frac{3}{2}}\right)\nonumber\\[5pt]
	&=k^{\frac{1}{2}}(k-1)^2\left(\frac{1}{3\cdot 8}+\frac{1}{3\cdot8^2k^{\frac{3}{2}}}+\frac{1}{12\cdot 8^3k^3}\right)\nonumber\\[5pt]
	&\quad+\gamma k^2(k-1)\left(\frac{1}{3}+\frac{1}{8k^{\frac{3}{2}}}+\frac{6-k^{-1}}{6\cdot8^2k^3}+\frac{2-k^{-1}}{6\cdot8^3k^{\frac{9}{2}}}\right)\nonumber\\[5pt]
	&\leq k^3\left(\frac{1}{24}+\frac{1}{192k^{\frac{3}{2}}}+\frac{1}{6144k^3}+\gamma \left(\frac{1}{3}+\frac{1}{8k^{\frac{3}{2}}}+\frac{1}{64r^3}+\frac{1}{1536k^{\frac{9}{2}}}\right)\right)\nonumber\\[5pt]
	&\leq k^3\left(\frac{1}{24}+\frac{1}{192\cdot4^{\frac{3}{2}}}+\frac{1}{6144\cdot 4^3}\right.\nonumber\\[5pt]
	&\left.\quad + 0.616\cdot\left(\frac{1}{3}+\frac{1}{8\cdot4^{\frac{3}{2}}}+\frac{1}{64\cdot4^3}+\frac{1}{1536\cdot 4^{\frac{9}{2}}}\right)\right)\quad (\text{by} \ k\geq 4)\nonumber\\[5pt]
	&< 0.26k^3:=c_k,
\end{align*}
and so
\begin{equation}\label{eq-I1}   \theta\sin\left(\mu\theta\right)\prod_{j=0}^{m}\prod_{l=1}^{k-1}\cos((jk+l)\theta)\geq \mu\theta^2\exp\left(-c_km^3\theta^2\right).
\end{equation}
Applying \eqref{eq-I1} to \eqref{defi-I1}, we deduce that when $k\geq 4$, $m\geq 8k^{\frac{3}{2}}$ and $0<\mu\leq\frac{k(k-1)(2m+1)}{2}$,
\begin{align}
	I_{k,m}^{(1)}(\mu)&=\int_0^{\frac{2\pi}{k(k-1)(2m+1)}}\theta\sin\left(\mu\theta\right)\prod_{j=0}^{m}\prod_{l=1}^{k-1}\cos((jk+l)\theta)\mathrm{d}\theta\nonumber\\[5pt]
	&\geq \int_0^{\frac{4}{k(k-1)(2m+1)}}\theta\sin\left(\mu\theta\right)\prod_{j=0}^{m}\prod_{l=1}^{k-1 }\cos((jk+l)\theta)\mathrm{d}\theta\nonumber\\[5pt]
	&\geq\int_0^{\frac{4}{k(k-1)(2m+1)}}\mu\theta^2\exp\left(-c_km^3\theta^2\right)\mathrm{d}\theta\nonumber\\[3pt]
	&=\left\{\int_0^{\infty}-\int_{\frac{4}{k(k-1)(2m+1)}}^{\infty}\right\}\mu\theta^2\exp\left(-c_km^3\theta^2\right)\mathrm{d}\theta\nonumber\\[5pt]
	&=\frac{\mu}{2c_k^{\frac{3}{2}}m^{\frac{9}{2}}}\left(\int_0^{\infty}v^{\frac{1}{2}}e^{-v}\mathrm{d}v-\int_{\frac{16c_km^3}{k^2(k-1)^2(2m+1)^2}}^{\infty}v^{\frac{1}{2}}e^{-v}\mathrm{d}v\right)\nonumber\\[5pt]
	&=\frac{\mu}{2c_k^{\frac{3}{2}}m^{\frac{9}{2}}}\left(\frac{\sqrt{\pi}}{2}-\int_{\frac{16c_km^3}{k^2(k-1)^2(2m+1)^2}}^{\infty}v^{\frac{1}{2}}e^{-v}\mathrm{d}v\right).\nonumber
\end{align}	
When $m\geq 8k^{\frac{3}{2}}$, we see that
\begin{align}
	\frac{16c_km^3}{k^2(k-1)^2(2m+1)^2}
	&\geq\frac{16\cdot 0.26k^3\cdot8^3k^\frac{9}{2}}{k^2(k-1)^2(2\cdot8k^{\frac{3}{2}}+1)^2}\nonumber\\[5pt]
	&\geq \frac{16\cdot 0.26k^3\cdot8^3k^\frac{9}{2}}{k^2k^2(17k^{\frac{3}{2}})^2}\nonumber\\[5pt]
	&= \frac{2129.92\sqrt{k}}{289}\nonumber\\[5pt]
	&\geq\frac{2129\sqrt{4}}{289}>14.7 \quad(\text{by} \ k\geq 4),
\end{align}
so
\[\int_{\frac{16c_km^3}{k^2(k-1)^2(2m+1)^2}}^{\infty}v^{\frac{1}{2}}e^{-v}\mathrm{d}v<\int_{14.7}^{\infty}v^{\frac{1}{2}}e^{-v}\mathrm{d}v<1.64\times 10^{-6}.\]
As a result, we can assert that   when $k\geq 4$, $m\geq 8k^{\frac{3}{2}}$ and $0<\mu\leq\frac{k(k-1)(2m+1)}{2}$,
\begin{align}
	I_{k,m}^{(1)}(\mu)> \frac{\mu}{2c_k^{\frac{3}{2}}m^{\frac{9}{2}}}\left(\frac{\sqrt{\pi}}{2}-1.64\times 10^{-6}\right)>\frac{3.34\mu}{k^{\frac{9}{2}}m^{\frac{9}{2}}}.
	\label{eq-I1-0}
\end{align}

We now turn to estimate the value of $I_{k,m}^{(2)}(\mu)$ given by
\begin{equation}\label{defi-i2}
	I_{k,m}^{(2)}(\mu)= \int_{\frac{2\pi}{k(k-1)(2m+1)}}^{\frac{\pi}{2km+2(k-1)}}\theta\sin\left(\mu\theta\right)\prod_{j=0}^{m}\prod_{l=1}^{k-1}\cos((jk+l)\theta)\mathrm{d}\theta.
\end{equation}
When $\frac{2\pi}{k(k-1)(2m+1)}\leq\theta\leq\frac{\pi}{2km+2(k-1)}$, we have
$0 \leq(j k+l)\theta\leq \frac{\pi}{2} \ \text{for}\ 0\leq j\leq m \ \text{and}\ 1\leq l\leq k-1.$
In light of \eqref{eq-2.7-1}, we derive that
\[\cos((jk+l)\theta)\leq \exp\left(-\frac{(jk+l)^2\theta^2}{2}\right).\]
Hence
\begin{align}
	&\left|\prod_{j=0}^{m}\prod_{l=1}^{k-1}\cos((jk+l)\theta)\right|\nonumber\\[3pt]
	&\quad\leq\exp\left(-\frac{1}{2}\theta^2\sum_{j=0}^{m}\sum_{l=1}^{k-1}(jk+l)^2\right)	\nonumber\\[5pt]
	&\quad=\exp\left(-\frac{1}{2}k(k-1)\theta^2\left(\frac{km^3}{3}+km^2+\frac{(6k-1)m}{6}+\frac{2k-1}{6}\right)\right)\nonumber\\[5pt] &\quad\leq\exp\left(-\frac{\pi^2}{2k(k-1)\left(m+\frac{1}{2}\right)^2}\left(\frac{km^3}{3}+km^2+\frac{(6k-1)m}{6}+\frac{2k-1}{6}\right)\right)\nonumber\\[5pt]
	& ~~~~~~~~~~~~~~~~~~~~~~~~~~~~~~~~~~~~~~~\left(\text{by}\ \frac{2\pi}{k(k-1)(2m+1)}\leq\theta\leq\frac{\pi}{2k m+2(k-1)}\right)
	\nonumber\\[5pt] &\quad=\exp\left(-\frac{\pi^2}{2k(k-1)}\cdot\frac{km}{3}\cdot\frac{m^2+3m+\frac{6k-1}{2k}+\frac{2k-1}{2km}}{m^2+m+\frac{1}{4}}\right)\nonumber\\[5pt]
	&\quad\leq\exp\left(-\frac{\pi^2}{2k(k-1)}\cdot\frac{km}{3}\right)=\exp\left(-\frac{\pi^2 m}{6(k-1)}\right)<\exp\left(-\frac{\pi^2 m}{6k}\right).
	\label {eq-i2-1}
\end{align}	
Applying  \eqref{eq-i2-1}  to \eqref{defi-i2}, and in view of \eqref{eq-2.6} and \eqref{eq-I1-0}, we derive that when $k\geq 4$, $m\geq 8k^{\frac{3}{2}}$ and $0<\mu\leq\frac{k(k-1)(2m+1)}{2}$,
\begin{align}\label {eq-i2}
	|I_{k,m}^{(2)}(\mu)|
	&\overset{\eqref{eq-2.6}}{\leq} \mu \exp\left(-\frac{\pi^2m}{6k}\right)\int_{\frac{2\pi}{k(k-1)(2m+1)}}^{\frac{\pi}{2km+2(k-1)}}\theta^2\mathrm{d}\theta	\nonumber\\[5pt]
	&\leq \frac{\mu\pi^3}{3}\left(\frac{1}{(2km+2(k-1))^3}-\frac{8}{(k(k-1)(2m+1))^3}\right)
	\exp\left(-\frac{\pi^2 m}{6k}\right)	\nonumber\\[5pt]
	&\leq\frac{\mu\pi^3}{3(2km+2(k-1))^3}\exp\left(-\frac{\pi^2 m}{6k}\right)	\nonumber\\[5pt]
	&\leq\frac{\mu\pi^3}{3(8m)^3}\exp\left(-\frac{\pi^2 m}{6k}\right) \quad \left(\text{by}\ k\geq
	4\right)\nonumber\\[5pt]
	&\overset{\eqref{eq-I1-0}}{\leq} \frac{\pi^{3}k^{\frac{9}{2}}m^{\frac{3}{2}}}{5130} \exp\left(-\frac{\pi^2 m}
	{6k}\right)I_{k,m}^{(1)}(\mu).
\end{align}	
Define
\[f_k(m):=\frac{\pi^{3}k^{\frac{9}{2}}m^{\frac{3}{2}}}{5130} \exp\left(-\frac{\pi^2 m}{6k}\right).\]
We claim that $f_k'(m)<0$ for $k\geq 4$ and  $m\geq 8k^{\frac{3}{2}}$. Since $f_k(m)>0$ for $k\geq 4$ and $m\geq 8k^{\frac{3}{2}}$, we have
\begin{equation}\label{eq-f'}
	\frac{\mathrm{d}}{\mathrm{d}m} f_k(m)=\frac{\mathrm{d}}{\mathrm{d}m} \mathrm{e}^{\ln{f_k(m)}}=f_k(m)\frac{\mathrm{d}}{\mathrm{d}m} \ln{f_k(m)}.
\end{equation}
Observe that when $k\geq 4$ and $m\geq 8k^{\frac{3}{2}}$,
\begin{align*}
	\frac{\mathrm{d}}{\mathrm{d}m} \ln{f_k(m)}&=\frac{3}{2m} -\frac{\pi^2 }{6k}\leq\frac{3}{2\cdot 8k^{\frac{3}{2}}} -\frac{\pi^2 }{6k}=\frac{\pi^2}{6k}\left(\frac{9}{8\pi^2k^{\frac{1}{2}}}- 1 \right)<0,
\end{align*}
and this yields that $f_k'(m)<0$ for $k\geq4$ and $m\geq 8k^{\frac{3}{2}}$ as claimed. Consequently,
\begin{align}\label{eq-f}
	f_k(m)\leq f_k(8k^{\frac{3}{2}})=\frac{8^{\frac{3}{2}}\pi^3}{5130}k^{\frac{27}{4}}\exp\left(-\frac{4\pi^2 k^{\frac{1}{2}}}{3}\right).
\end{align}
Applying  \eqref{eq-f} to \eqref{eq-i2}, we obtain
\begin{align}\label {eq-i2-2}
	|I_{k,m}^{(2)}(\mu)|
	\leq\frac{8^{\frac{3}{2}}\pi^3}{5130}k^{\frac{27}{4}}\exp\left(-\frac{4\pi^2 k^{\frac{1}{2}}}{3}\right)I_{k,m}^{(1)}(\mu).
\end{align}	
Define
\[h_1(k):=\exp\left(-\frac{4\pi^2 k^{\frac{1}{2}}}{3}\right)k^{\frac{27}{4}}.\]
Since $h_1(k)>0$ for $k\geq 4$, we find that
\begin{equation}\label{eq-h_1}
	\frac{\mathrm{d}}{\mathrm{d}k} {h_1(k)}=\frac{\mathrm{d}}{\mathrm{d}k} e^{\ln h_1(k)}=h_1(k)\frac{\mathrm{d}}{\mathrm{d}k} \ln{h_1(k)},
\end{equation}
and since
\begin{align}
	\frac{\mathrm{d}}{\mathrm{d}k} \ln{h_1(k)}
	&=\frac{27}{4k}-\frac{2\pi^2}{3k^{\frac{1}{2}}} \nonumber\\[5pt]
	&=\frac{1}{k}\left(\frac{27}{4}-\frac{2\pi^2k^{\frac{1}{2}}}{3}\right)\nonumber\\[5pt]
	&\leq \frac{1}{k}\left(\frac{27}{4}-\frac{4\pi^2}{3}\right)\quad(\text{by}\ k\geq 4)\nonumber\\[5pt]
	&<-\frac{6}{k}<0,\nonumber
\end{align}
it follows that $h_1'(k)<0$ for $k\geq 4$. Hence  $h_1(k)\leq h_1(4)$ for $k\geq 4$. Therefore,
\begin{align}\label{eq-i2-3}
	|I_{k,m}^{(2)}(\mu)|&\overset{\eqref{eq-i2-2}}{\leq}\frac{ 8^{\frac{3}{2}}\pi^3}{5130}\exp\left(-\frac{8\pi^2}{3} \right)\cdot4^{\frac{27}{4}}I_{k,m}^{(1)}(\mu)\nonumber\\[5pt]
	&<5.89\times10^{-9} I_{k,m}^{(1)}(\mu).
\end{align}

Finally, we turn to  estimate the value of $I_{k,m}^{(3)}(\mu)$ defined by
\begin{equation}\label{defi-i3}
	I_{k,m}^{(3)}(\mu)
	=\int_{\frac{\pi}{2km+2(k-1)}}^{\frac{\pi}{2}}\theta\sin\left(\mu\theta\right)\prod_{j=0}^{m}\prod_{l=1}^{k-1}\cos((j k+l)\theta)\mathrm{d}\theta.
\end{equation}
Let $C=\left\{\frac{i\pi}{2k} \vert i=1,2,\ldots,k\right\}$, it is easy to see that
\[\int_C\theta\sin\left(\mu\theta\right)\prod_{j=0}^{m}\prod_{l=1}^{k-1}\cos((j k+l)\theta)\mathrm{d}\theta=0,\]
so
\begin{equation}\label{defi-i3-1}
	I_{k,m}^{(3)}(\mu)
	=\int_{\left[\frac{\pi}{2km+2(k-1)},\frac{\pi}{2}\right] \setminus C}\theta\sin\left(\mu\theta\right)\prod_{j=0}^{m}\prod_{l=1}^{k-1}\cos((j k+l)\theta)\mathrm{d}\theta.
\end{equation}
When  $\frac{\pi}{2km+2(k-1)}\leq \theta \leq \frac{\pi}{2}$ and $\theta\neq \frac{i\pi}{2k}$ ($i=1,2,\ldots,k$),
by  \eqref{eq-2.7}, \eqref{eq-2.8} and \eqref{eq-2.9}, we deduce that
\begin{align}
	&\left|\prod_{j=0}^{m}\prod_{l=1}^{k-1}\cos((jk+l)\theta)\right|\nonumber\\[5pt] &\quad\overset{\eqref{eq-2.7}}{\leq}\exp\left(-\frac{1}{2}\sum_{j=0}^{m}\sum_{l=1}^{k-1}\sin^2((jk+l)\theta)-\frac{1}{4}\sum_{j=0}^{m}\sum_{l=1}^{k-1}\sin^4((jk+l)\theta)\right)	\nonumber\\[5pt]
	&\quad=\exp\left (-\frac{1}{2}\left(\sum_{j=1}^{km+k-1}\sin^2(j\theta)-\sum_{j=1}^{m}\sin^2(jk\theta)\right)\right.\nonumber\\[5pt]
	&\left.\quad\quad-\frac{1}{4}\left(\sum_{j=1}^{km+k-1}\sin^4(j\theta)-\sum_{j=1}^{m}\sin^4(jk\theta)\right)\right )\nonumber\\[5pt] &\quad\overset{\eqref{eq-2.8}\&\eqref{eq-2.9}}{=}\exp\left(-\frac{11(k-1)(m+1)}{32}+\frac{3\sin((2km+2k-1)\theta)}{16\sin(\theta)}\right.\nonumber\\[5pt]
	&\left.\quad\quad -\frac{\sin((2km+2k-1)2\theta)}{64\sin(2\theta)}-\frac{3\sin((2m+1)k\theta)}{16\sin(k\theta)}+\frac{\sin((2m+1)2k\theta)}{64\sin(2k\theta)}\right)\nonumber\\[5pt]
	&\quad:=E_{k,m}(\theta).\label{defi-e}
\end{align}	
We claim that for  $k\geq 4 $,  $ m\geq 8k^{\frac{3}{2}}$ and $\frac{\pi}{2km+2(k-1)}\leq \theta \leq \frac{\pi}{2}$ (where $\theta\neq \frac{i\pi}{2k}$, $i=1,2,\ldots,k$),
\begin{equation}\label{eval-value-e}
	E_{k,m}(\theta)<\exp\left(-0.381m-0.224\right) .
\end{equation}
We approach the proof of  \eqref{eval-value-e} through a two-step process. First, we consider the interval $\frac{\pi}{2km+2(k-1)}\leq \theta < \frac{\pi}{2k}$. Since
$\frac{\pi}{2km+2(k-1)}\leq\theta<2\theta<k\theta<\frac{\pi}{2}$,  by \eqref{eq-2.6}, we get that,
\begin{align}\label{eq-2.37}
	\sin(i\theta)&\geq\sin\left(\frac{i\pi}{2km+2(k-1)}\right)\nonumber\\[5pt]
	&\geq\frac{i\pi}{2km+2(k-1)}-\frac{\left(\frac{i\pi}{2km+2(k-1)}\right)^3}{6}\nonumber\\[5pt]
	&\geq \frac{i\pi}{2km+2(k-1)}\left(1-\frac{\left(\frac{k\pi}{2km+2(k-1)}\right)^2}{6}\right),
\end{align}
where $i=1,\,2,\,k$.
Applying \eqref{eq-2.7-2} and \eqref{eq-2.37} in \eqref{defi-e},  we obtain
\begin{align}
	E_{k,m}(\theta)
	&\leq\exp\left (-\frac{11(k-1)(m+1)}{32}+\frac{3}{16\sin(\theta)}+\frac{1}{64\sin(2\theta)}+\frac{3}{16\sin(k\theta)}\right.\nonumber\\[5pt]
	&\left.\quad+\left|\frac{\sin((2m+1)2k\theta)}{64\sin(2k\theta)}\right|\right )\nonumber\\[5pt]	
	&\overset{\eqref{eq-2.37}\&\eqref{eq-2.7-2}}{\leq}\exp\left (-\frac{11(k-1)(m+1)}{32}+\frac{2m+1}{64}+\frac{3}{16\left(\frac{\pi}{2km+2(k-1)}\left(1-\frac{\left(\frac{k\pi}{2km+2(k-1)}\right)^2}{6}\right)\right)}\right.\nonumber\\[5pt]	&\left.\quad+\frac{1}{64\left(\frac{2\pi}{2km+2(k-1)}\left(1-\frac{\left(\frac{k\pi}{2km+2(k-1)}\right)^2}{6}\right)\right)}+\frac{3}{16\left(\frac{k\pi}{2km+2(k-1)}\left(1-\frac{\left(\frac{k\pi}{2km+2(k-1)}\right)^2}{6}\right)\right)}\right)\nonumber\\[5pt]
	&=\exp\left (\frac{(12-11k)m}{32}+\frac{23-22k}{64}+\frac{24+25k}{128k\left(\frac{\pi}{2km+2(k-1)}\left(1-\frac{\left(\frac{k\pi}{2km+2(k-1)}\right)^2}{6}\right)\right)}\right )\nonumber\\[5pt]
	&=\exp\left (\frac{(12-11k)m}{32}+\frac{23-22k}{64}+\frac{(24+25k)(2km+2(k-1))}{128\pi k\left(1-\frac{\pi^2k^2}{6(2km+2(k-1))^2}\right)}\right ).\nonumber
\end{align}	
When  $k\geq 4$ and $m\geq 8k^{\frac{3}{2}}$, we have
\begin{align}
	1-\frac{\pi^2k^2}{6(2km+2(k-1))^2}
	&\geq1-\frac{\pi^2k^2}{6\left(16k^{\frac{5}{2}}+2(k-1)\right)^2}\quad(\text{by}\ m\geq 8k^{\frac{3}{2}})\nonumber\\[5pt]
	&=1-\frac{\pi^2}{6\left(16k^{\frac{3}{2}}+2-2k^{-1}\right)^2}\nonumber\\[5pt]
	&\geq 1-\frac{\pi^2}{6\left(16\cdot4^{\frac{3}{2}}+2-\frac{1}{2}\right)^2}\quad(\text{by}\ k\geq 4)\nonumber\\[5pt]
	&=1-\frac{\pi ^2}{100621.5}
	>0.9999.\nonumber
\end{align}
It follows that for  $k\geq 4$,  $m\geq 8k^{\frac{3}{2}}$ and $\frac{\pi}{2km+2(k-1)}\leq \theta < \frac{\pi}{2k}$,
\begin{align}\label{eq-E-1}
	E_{k,m}(\theta)&\leq\exp\left (\frac{(12-11k)m}{32}+\frac{23-22k}{64}+\frac{(24+25k)(2km+2(k-1))}{0.9999\cdot128\pi k}\right )\nonumber\\[5pt]
	&=\exp\left(\left(\frac{12-11k}{32}+\frac{24+25k}{0.9999\cdot64\pi}\right)m+\frac{23-22k}{64}+\frac{(24+25k)(1-k^{-1})}{0.9999\cdot64\pi}\right)\nonumber\\[5pt]
	&\leq\exp\left((0.495-0.219k)m+0.479-0.219k\right) \nonumber\\[5pt]
	&\leq\exp\left(-0.381m-0.397\right)\quad (\text{by}\ k\geq 4).
\end{align}

Next we consider the interval  $\frac{\pi}{2k}\leq \theta \leq \frac{\pi}{2}$ and $\theta\neq \frac{i\pi}{2k}$ ($i=1,2,\ldots,k$). Employing  \eqref{eq-2.6} and \eqref{eq-2.7-2}, we deduce that
\begin{align}\label{eq-E-2}
	E_{k,m}(\theta)&\leq\exp\left(-\frac{11(k-1)(m+1)}{32}+\frac{3}{16\sin\left(\theta\right)}\right.\nonumber\\[5pt] &\left.\quad\quad+\left|\frac{\sin((2km+2k-1)2\theta)}{64\sin(2\theta)}\right|+\left|\frac{3\sin((2m+1)k\theta)}{16\sin(k\theta)}\right|+\left|\frac{\sin((2m+1)2k\theta)}{64\sin(2k\theta)}\right|\right)\nonumber\\[5pt]
	&\overset{\eqref{eq-2.7-2}}{\leq}\exp\left(-\frac{11(k-1)(m+1)}{32}+\frac{3}{16\sin\left(\frac{\pi}{2k}\right)}\right.\nonumber\\[5pt]
	&\left.\quad\quad+\frac{2km+2k-1}{64}+\frac{3(2m+1)}{16}+\frac{2m+1}{64}\right)\nonumber\\[5pt] &\overset{\eqref{eq-2.6}}{\leq}\exp\left(-\frac{11(k-1)(m+1)}{32}+\frac{3}{16\left(\frac{\pi}{2k}\left(1-\frac{\left(\frac{\pi}{2k}\right)^2}{6}\right)\right)}\right.\nonumber\\[5pt]
	&\left.\quad\quad+\frac{2km+2k-1}{64}+\frac{3(2m+1)}{16}+\frac{2m+1}{64}\right)\nonumber\\[5pt] &=\exp\left(\left(\frac{3}{4}-\frac{5k}{16}\right)m-\frac{5k}{16}+\frac{17}{32}+\frac{3k}{8\pi\left(1-\frac{\pi^2}{24k^2}\right)}\right)\nonumber\\[5pt]
	&\leq\exp\left(\left(\frac{3}{4}-\frac{5k}{16}\right)m-\frac{5 k}{16}+\frac{17}{32}+\frac{3k}{0.9742\cdot8\pi}\right)\quad(\text{by}\ k\geq 4)\nonumber\\[5pt]
	&\leq\exp\left(\left(\frac{3}{4}-\frac{5k}{16}\right)m+\frac{17}{32}-0.189k\right)\nonumber\\[5pt]
	&\leq\exp\left(-0.5m-0.224\right)\quad(\text{by}\ k\geq 4).
\end{align}	

Combining \eqref{eq-E-1} and \eqref{eq-E-2} yields \eqref{eval-value-e}, so the claim is verified.
Substituting \eqref{eval-value-e} to \eqref{defi-i3-1}, and in view of \eqref{eq-2.6} and \eqref{eq-I1-0}, we derive that
\begin{align}\label{eq-3}
	|I_{k,m}^{(3)}(\mu)|
	&\overset{\eqref{eq-2.6}}{\leq }\mu \exp\left(-0.381m-0.224\right)\int_{\frac{\pi}{2km+2(k-1)}}^{\frac{\pi}{2}}\theta^2\mathrm{d}\theta	\nonumber\\[5pt]
	&\leq \frac{\mu\pi^3}{24} \exp\left(-0.381m-0.224\right)	\nonumber\\[5pt]
	&\overset{\eqref{eq-I1-0}}{\leq} \frac{\pi^3k^{\frac{9}{2}}m^{\frac{9}{2}}}{3.34\cdot 24} \exp\left(-0.381m-0.224\right)I_{k,m}^{(1)}(\mu).
\end{align}	
Define
\[g_k(m):=\frac{\pi^3k^{\frac{9}{2}}m^{\frac{9}{2}}}{3.34\cdot 24} \exp\left(-0.381m-0.224\right).\]
Since when $k\geq 4$ and $m\geq 8k^{\frac{3}{2}}$,  we have $g_k(m)>0$ and
\begin{align}
	\frac{\mathrm{d}}{\mathrm{d}m} {g_k(m)}&=\frac{\mathrm{d}}{\mathrm{d}m} e^{\ln g_k(m)} \nonumber\\[5pt]
	&=g_k(m)\frac{\mathrm{d}}{\mathrm{d}m} \ln{g_k(m)}\nonumber\\[5pt]
	& =g_k(m)\left(\frac{9}{2m}-0.381\right)\nonumber\\[5pt]
	&\leq g_k(m)\left(\frac{9}{2\cdot8\cdot 4^{\frac{3}{2}}}-0.381\right)\nonumber\\[5pt]
	&<-0.31g_k(m)<0,\nonumber
\end{align}
it follows that $g_k'(m)<0$ when $k\geq 4$ and $m\geq 8k^{\frac{3}{2}}$, and so for $k\geq 4$,
\begin{align}
	g_k(m)&\leq g_k(8k^{\frac{3}{2}})=\frac{8^{\frac{9}{2}}\pi^3k^{\frac{45}{4}}}{3.34\cdot 24} \exp\left(-3.048k^{\frac{3}{2}}-0.224\right).
	\label{eq-i-3}
\end{align}
Define
\[h_2(k):=\exp\left(-3.048k^{\frac{3}{2}}-0.224\right)k^{\frac{45}{4}}.\]
When $k\geq 4$, we have  $h_2(k)>0$ and
\begin{align}
	\frac{\mathrm{d}}{\mathrm{d}k} {h_2(k)}
	&=\frac{\mathrm{d}}{\mathrm{d}k} e^{\ln h_2(k)}\nonumber\\[5pt]
	& =h_2(k)\frac{\mathrm{d}}{\mathrm{d}k} \ln{h_2(k)}\nonumber \\[5pt]
	&=h_2(k)\left(\frac{45}{4k}-3.048\cdot\frac{3k^{\frac{1}{2}}}{2}\right)\nonumber\\[5pt]
	&\leq h_2(k)\left(\frac{45}{4\cdot 4}-3.048\cdot\frac{3\sqrt{4}}{2}\right)\quad(\text{by} \ k\geq 4) \nonumber\\[5pt]
	&<-6.3h_2(k)<0,\nonumber
\end{align}
so $h_2'(k)<0$ for $k\geq 4$, and hence  for $k\geq 4$,
\begin{align}
	g_k(m)\leq\frac{8^{\frac{9}{2}}\pi^3}{3.34\cdot 24}\exp\left(-3.048\cdot4^{\frac{3}{2}}-0.224\right)\cdot4^{\frac{45}{4}}
	< 0.55.  \label{eval-val}
\end{align}
Substituting \eqref{eval-val} into \eqref{eval-integral}, we have
\begin{equation}\label{eq-i3-2}
	|I_{k,m}^{(3)}(\mu)|<0.55 I_{k,m}^{(1)}(\mu).
\end{equation}

Combining \eqref{eq-i2-3} and \eqref{eq-i3-2} yields \eqref{eq-4.5}, and so \eqref{eq-4.4} is valid. This leads to \eqref{eq-4.2} holds for $k\geq 4$, $m\geq 8k^{\frac{3}{2}}$ and  $\frac{k(k-1)m^2}{4}\leq n\leq\frac{k(k-1)(m+1)^2}{4}$, and so Lemma \ref{main-lemma} is verified. \qed
\section{Proofs of Theorem \ref{thm-k4} and Theorem \ref{main-thm}}

This section is devoted to the proofs of Theorem \ref{thm-k4} and Theorem \ref{main-thm}.  Prior to that, we demonstrate  the symmetry of $D_{k,m}(q)$.

\begin{thm}\label{thm-sys} For $k\geq 0$, the polynomials $D_{k,m}(q)$ are symmetric.
\end{thm}

\pf  Replacing $q$ by $q^{-1}$ in \eqref{eq-main}, we find that
\begin{align*}
	D_{k,m}(q^{-1})&=\prod_{j=0}^m\left(1+q^{-(kj+1)}\right)\left(1+q^{-(kj+2)}\right)\cdots\left(1+q^{-(kj+k-1)}\right)\nonumber\\[5pt]
	&=q^{-N(k,m)}\prod_{j=0}^m \left(1+q^{kj+1}\right)\left(1+q^{kj+2}\right)\cdots\left(1+q^{kj+k-1}\right) \\[5pt]
	&=q^{-N(k,m)} D_{k,m}(q)
\end{align*}
To wit,
\[D_{k,m}(q)=q^{N(k,m)}D_{k,m}(q^{-1}),\]
from which,  it follows that $D_{k,m}(q)$ is  symmetric. This completes the proof. \qed

We  give an inductive  proof of   Theorem \ref{thm-k4}   with the aid of Lemma \ref{main-lemma}.

\noindent{\it Proof of Theorem \ref{thm-k4}:} From Theorem \ref{thm-sys}, we see that   $D_{4,m}(q)$ is symmetric. Hence in order to prove Theorem \ref{thm-k4}, it suffices to show that
\begin{equation}\label{eq-a4-main}
	d_{4,m}(n)\geq d_{4,m}(n-1) \quad \text{for}\ m\geq 0,\  1\leq n\leq 3(m+1)^2 \ \text{and} \ n\not=4.
\end{equation}
Recall that $d_{4,m}(n)$ counts the number of $4$-regular partition into distinct parts where   the largest part is at most $4m+3$,  it is easy to check that for $m\geq 0$,
\begin{equation}\label{value-a-4}
	d_{4,m}(0)=d_{4,m}(1)=d_{4,m}(2)=1,\ d_{4,m}(3)=2,\ d_{4,m}(4)=1.
\end{equation}
Here we assume that $d_{4,m}(n)=0$ when $n<0$.  It  can be checked that \eqref{eq-a4-main} holds when $0\leq m\leq 63$.
In the following, we will demonstrate its validity for the case when $m\geq 64$. However,  our main objective is to show that  when $m\geq 64$,
\begin{equation}\label{eq-a4-maina}
	d_{4,m}(n)\geq   d_{4,m}(n-1), \quad 5\leq n\ \leq 12m+20
\end{equation}
and
\begin{equation}\label{eq-a4-mainb}
	d_{4,m}(n)\geq  d_{4,m}(n-1)+1, \quad 12m+21\leq n\ \leq 3(m+1)^2.
\end{equation}
which are immediate led to \eqref{eq-a4-main}.
It can be checked that \eqref{eq-a4-maina} and \eqref{eq-a4-mainb} are valid when $m=64$.
It remains to show that \eqref{eq-a4-maina} and \eqref{eq-a4-mainb} hold when $m > 64$. We  proceed by induction on $m$.
Assume that \eqref{eq-a4-maina} and \eqref{eq-a4-mainb} are valid for $m-1$, namely
\begin{equation}\label{eq-a4-mainat}
	d_{4,m-1}(n)\geq   d_{4,m-1}(n-1), \quad 5\leq n\ \leq 12m+8
\end{equation}
and
\begin{equation}\label{eq-a4-mainbt}
	d_{4,m-1}(n)\geq  d_{4,m-1}(n-1)+1, \quad 12m+9\leq n\ \leq 3m^2.
\end{equation}
We aim to show that \eqref{eq-a4-maina} and \eqref{eq-a4-mainb} hold.

Comparing coefficients of $q^n$ in
\begin{equation*}
	D_{4,m}(q)=\left(1+q^{4m+1}\right)\left(1+q^{4m+2}\right)\left(1+q^{4m+3}\right)D_{4,m-1}(q),
\end{equation*}
we obtain the following recurrence relation:
\begin{align}\label{equ-rec-a4}
	d_{4,m}(n)
	&=d_{4,m-1}(n)+d_{4,m-1}(n-4m-1)+d_{4,n-1}(n-4m-2)\nonumber\\[5pt]
	&\quad+d_{4,m-1}(n-4m-3)+d_{4,m-1}(n-8m-3)
	+d_{4,m-1}(n-8m-4)\nonumber\\[5pt]
	&\quad+d_{4,m-1}(n-8m-5)+d_{4,m-1}(n-12m-6),
\end{align}
thereby leading to
\begin{align}\label{equ-m-m-1}
	d_{4,m}(n)- d_{4,m}(n-1)&=d_{4,m-1}(n)- d_{4,m-1}(n-1)\nonumber\\[5pt]
	&\quad+d_{4,m-1}(n-4m-1)- d_{4,m-1}(n-4m-4)\nonumber\\[5pt]
	&\quad+d_{4,m-1}(n-8m-3)- d_{4,m-1}(n-8m-6)\nonumber\\[5pt]
	&\quad+d_{4,m-1}(n-12m-6)- d_{4,m-1}(n-12m-7).
\end{align}
When $5\leq n\ \leq 12m+20$ and $n\neq 12m+10$, applying  \eqref{eq-a4-mainat} and \eqref{eq-a4-mainbt}  to \eqref{equ-m-m-1},     we see  that
\[d_{4,m}(n)- d_{4,m}(n-1)\geq 0.\]
When $n=12m+10$, we observe that
\[d_{4,m-1}(n-12m-6)- d_{4,m-1}(n-12m-7)=d_{4,m-1}(4)- d_{4,m-1}(3)=-1.\]
But by \eqref{eq-a4-mainbt}, we have
\[d_{4,m-1}(n)- d_{4,m-1}(n-1)=d_{4,m-1}(12m+10)- d_{4,m-1}(12m+9)\geq 1,\]
which leads to  $ d_{4,m}(n)- d_{4,m}(n-1)\geq 0$ when $n=12m+10$. To sum up, we get
\[ d_{4,m}(n)- d_{4,m}(n-1)\geq 0, \quad  5\leq n\ \leq12m+20,\]
and so \eqref{eq-a4-maina} is valid.
Applying \eqref{eq-a4-mainat} and \eqref{eq-a4-mainbt}  to \eqref{equ-m-m-1} again, we  infer that
\begin{equation}\label{eq-a4-maintt}
	d_{4,m}(n)- d_{4,m}(n-1)\geq 1, \quad  12m+21\leq n\ \leq 3m^2.
\end{equation}
In view of Lemma \ref{main-lemma}, we see that
\begin{equation}\label{eq-a4-maintt2}
	d_{4,m}(n)- d_{4,m}(n-1)\geq 1, \quad 3m^2< n\leq 3(m+1)^2.
\end{equation}
Combining \eqref{eq-a4-maintt} and \eqref{eq-a4-maintt2}, we confirm that \eqref{eq-a4-mainb} holds. Together with \eqref{eq-a4-maina}, we deduce \eqref{eq-a4-main} holds, and so $D_{4,m}(q)$ is unimodal, except at the coefficients of $q^4$ and $q^{N(4,m)-4}$. This completes the proof of Theorem \ref{thm-k4}.\qed

We conclude this paper with the proof of Theorem \ref{main-thm}  by the utilization of Lemma \ref{main-lemma}.

\noindent{\it Proof of Theorem \ref{main-thm}:}
Given $k\geq5$ and $m_0\geq 0$, assume that $D_{k,m_0}(q)$ is unimodal. We proceed to  show that the polynomial $D_{k,m}(q)$ is unimodal for $m\geq m_0$  by induction on $m$. Considering the symmetry of $D_{k,m}(q)$, it suffices to show that  for $m> m_0$ and $1\leq n\leq   \lfloor \frac{k(k-1)(m+1)^2}{4}\rfloor$,
\begin{equation}\label{ineq-n-k-2}
	d_{k,m}(n)\geq d_{k,m}(n-1).
\end{equation}
Assume that  \eqref{ineq-n-k-2} is valid for $m-1$, that is,   for  $m>  m_0$ and $1\leq n\leq \lfloor \frac{k(k-1)m^2}{4}\rfloor$,
\begin{equation}\label{ineq-n-k-1}
	d_{k,m-1}(n)\geq d_{k,m-1}(n-1).
\end{equation}
We intend to show that \eqref{ineq-n-k-2} holds for $m> m_0$ and $1\leq n\leq \lfloor \frac{k(k-1)(m+1)^2}{4}\rfloor$.
By comparing the coefficients of $q^n$ in the polynomial
\begin{align*}
	D_{k,m}(q)= \left(1+q^{km+1}\right)\left(1+q^{km+2}\right)\cdots\left(1+q^{km+k-1}\right)D_{k,m-1}(q),
\end{align*}
it can be determined that
\begin{align*}
	d_{k,m}(n)= \sum_{\substack{i_j=0 \ \text{or}\ km+j\\1\leq j\leq k-1}}d_{k,m-1}(n-i_1-\cdots-i_{k-1}),
\end{align*}
which  leads to
\begin{align}\label{equ-rec-k}
	&d_{k,m}(n)-d_{k,m}(n-1)\nonumber\\[5pt]
	&\quad=\sum_{\substack{i_j=0 \ \text{or}\ km+j\\1\leq j\leq k-1}}\left(d_{k,m-1}(n-i_1-\cdots-i_j)-d_{k,m-1}(n-i_1-\cdots-i_j-1)\right).
\end{align}
Utilizing \eqref{ineq-n-k-1} in \eqref{equ-rec-k} yields that the validity of \eqref{ineq-n-k-2} for $m>m_0$ and $1\leq n\leq \lfloor \frac{k(k-1)m^2}{4}\rfloor$.
In view of Lemma \ref{main-lemma}, we see that \eqref{ineq-n-k-2}  holds for $m\geq 8k^{\frac{3}{2}}$ and $\lceil \frac{k(k-1)m^2}{4}\rceil\leq n\leq \lfloor \frac{k(k-1)(m+1)^2}{4}\rfloor$. Given the condition that  \eqref{ineq-n-k-2} holds for $m_0< m< 8k^{\frac{3}{2}}$ and  $\lceil \frac{k(k-1)m^2}{4}\rceil\leq n\leq \lfloor \frac{k(k-1)(m+1)^2}{4}\rfloor$, we reach the conclusion that \eqref{ineq-n-k-2} is valid for $m> m_0$ and $1\leq n\leq \lfloor \frac{k(k-1)(m+1)^2}{4}\rfloor$. Therefore,  $D_{k,m}(q)$ is unimodal for  $m\geq m_0$. Thus, we complete the proof of Theorem \ref{main-thm}.
\qed

\vskip 0.2cm
		\noindent{\bf Acknowledgment.} This work
		was supported by   the National Science Foundation of China.

\end{document}